\documentclass{amsproc}
\usepackage{amsmath,amsfonts}

\newtheorem{thm}{Theorem}[section]

\begin{document}

\def\l{\lambda}
\def\t{\theta}
\def\T{\Theta}
\def\m{\mu}
\def\a{\alpha}
\def\b{\beta}
\def\g{\gamma}
\def\o{\omega}
\def\p{\varphi}

\def\psi{{\rho}}

\def\O{\Omega}

\def\N{{\mathbb N}}
\def\Z{{\mathbb Z}}

\def\L{{\mathcal L}}

\def\R{{\mathbb R}}
\def\P{{\mathbb P}}

\title{The $M/M/1$ queue is Bernoulli}

\author{Michael Keane}
\address{Department of Mathematics and Computer Science, Wesleyan University, Middletown, CT 06459, USA}
\author{Neil O'Connell}
\address{Mathematics Institute, University of Warwick, Coventry CV4 7AL, UK} 

\maketitle

\begin{abstract}
The classical output theorem for the $M/M/1$ queue, due to Burke (1956), 
states that the departure process from a stationary $M/M/1$ queue, in equilibrium, 
has the same law as the arrivals process, that is, it is a Poisson process.
In this paper we show that the associated measure-preserving transformation
is metrically isomorphic to a two-sided Bernoulli shift.  We also discuss some extensions
of Burke's theorem where it remains an open problem to determine if, or under
what conditions, the analogue of this result holds.

\end{abstract}


\section{Introduction}

The classical output theorem for the $M/M/1$ queue,
due to Burke~\cite{burke}, states that
the departure process from a stationary $M/M/1$ queue, in equilibrium, has
the same law as the arrivals process, that is, it is a Poisson process.
To be more precise, let $A$ and $S$ be Poisson processes on $\R$
with repective intensities $\lambda<\xi$ and define, for $t\in\R$,
$$Q(t)=\sup_{s\le t} (A(s,t]-S(s,t]).$$
For each $t$, $Q(t)$ should be interpreted as the number of customers
in the queue at time $t$.  Customers arrive according to the Poisson
process $A$ (the arrivals process) and at times given by the points 
of $S$, if the queue is non-empty, a customer is served and departs
from the queue.  The departure process $D$ is defined to be the point
process of times at which customers depart from the queue.  More
precisely, we define, for $s<t$,
$$D(s,t]=Q(s)+A(s,t]-Q(t).$$
Burke's theorem states that $D$ is a Poisson process with intensity $\lambda$,
and moreover that $(D(t,0],t<0)$ is independent of $Q(0)$.
The standard proof of this fact, due to Reich~\cite{reich}, is a reversibility argument 
which exploits the dynamical symmetry of the queue and the fact that $Q$ is a 
stationary, reversible Markov chain.  For more background on queueing theory, 
see, for example, Kelly~\cite{kelly}.

The nature of Burke's theorem suggest that there may be a measure-preserving
transformation somewhere nearby. It is not immediately obvious how to
find it, since $D$ is not only a function of $A$, it also depends on $S$.  
However, it was shown in~\cite{oy02}
that, if we define $R=A+S-D$, then the pair $(D,R)$
has the same joint law as $(A,S)$, thus exhibiting a measure-preserving
transformation; moreover, the restriction of $(D,R)$ to $(-\infty,0]^2$ is independent
of $Q(0)$.  We can restate this as follows.  For $t\in\R$, set
$$X(t)=\begin{cases} S(0,t]-A(0,t] & t>0\\
A(-t,0]-S(-t,0] & t\le 0.\end{cases}$$
and
$$Y(t)=\begin{cases} R(0,t]-D(0,t] & t>0\\
D(-t,0]-R(-t,0] & t\le 0.\end{cases}$$
Note that we can write 
$$Y(t)=2M(t)-X(t)-2M(0),\qquad M(t)=\sup_{-\infty<s\le t} X(s).$$
Then $X$ is a two-sided continuous-time simple random walk
with positive drift $\xi-\lambda$, and the transformation
which maps $X$ to $Y$ is measure-preserving; moreover, $(Y(t),t\le 0)$ is
independent of $Q(0)\equiv M(0)$. This statement can be further simplified 
by considering only the times at which events occur (i.e.
the times at which the random walk $X$ jumps).  Denote these
times (which are almost surely distinct) by 
$$\cdots<\tau_{-2}<\tau_{-1}<0<\tau_1<\tau_2<\cdots$$
and set $x_n=X(\tau_n)$ and $y_n=Y(\tau_n)$, for $n\in\Z$.
Note that, for $n\in\Z$, 
$$y_n=2s_n-x_n-2s_0,\qquad s_n=\sup_{m\le n}x_m.$$
Then $(x_n,n\in\Z)$ is a two-sided, discrete-time simple random
walk, as is $(y_n,n\in\Z)$, and $(y_n,n\le 0)$ is independent of $s_0$.
Finally, let $\O=\{-1,1\}^\Z$ be equipped with Bernoulli product measure with 
parameter $p=\xi/(\lambda+\xi)$.  Set $\epsilon_n=x_n-x_{n-1}$ and 
$\sigma_n=y_n-y_{n-1}$.  Then we can write $\sigma=T\epsilon$, where $T$, 
defined almost everywhere on $\O$, is a measure-preserving transformation.
The fact that $(\sigma_n,n\le 0)$ is independent of $s_0$ can now be interpreted
as saying that $T$ has a {\em factor} which is Bernoulli, that is, a factor which
is metrically isomorphic to a two-sided Bernoulli shift (see section 2 for details).  
The main result of this paper is that $T$ is, in fact, Bernoulli.  
This will be presented in section 2.
In section 3 we discuss the Brownian analogue of Burke's theorem where it is only
possible to show that the corresponding transformation has a Bernoulli factor.
The difficulty here is similar to that encountered in the open question,
posed by Marc Yor, of determining whether L\'evy's transformation of 
Brownian motion is ergodic.  Dubins and Smorodinsky~\cite{ds} proved
that there is a discrete version of L\'evy's transformation which is
isomorphic to a one-sided Bernoulli shift.
In section 4 we describe a natural extension of Burke's theorem
to the more general setting of iterated random functions,
and leave it as an open problem to determine under what conditions
the corresponding transformation is Bernoulli.

\section{The main result}

Let $\m$ be a Bernoulli product measure on $\O=\{-1,1\}^\Z$ with
$$\m\{\o\in\O:\ \o_0=1\}=p>1/2.$$  Define a two-sided 
simple random walk $x=(x_n,n\in\Z)$ by $x_0=0$,
$$x_n=\begin{cases}x_{n-1}+\o_n & n>0 ,\\
x_{n+1}-\o_{n+1} & n<0.\end{cases}$$
For $n\in\Z$, set $s_n=\sup_{m\le n} x_m$ and $\O'=\{s_0(\o)<\infty\}$.
Note that $\m(\O')=1$.  Write $y=2s-x$ and define $T:\O'\to\O$ by setting
$(T\o)_n=y_n-y_{n-1}$ for each $n\in\Z$.  In order to discuss the inverse
transformation we further define $$\O''=\{\o\in\O':\ \liminf_n(s_n-x_n)=0\}$$
and note that $\m(\O'')=1$.  Let $R:\O\to\O$ be the `time-reversal' operator
defined by $(R\o)_n=\o_{-n}$ for $n\in\Z$, and set $\rho=p^{-1}(1-p)$.
We first recall the analogue of Burke's theorem in this discrete setting.
\begin{thm}\label{burke}
\ \
\begin{itemize}
\item[i.]  $\m\circ T^{-1}=\mu$.
\item[ii.] For $x\ge 0$, $\m\{\o:\ s_0(\o)=x\}=(1-\rho)\rho^x$.
\item[iii.] The random variable $s_0$ is independent of $((T\o)_n,n\le 0)$.
\item[iv.] If $\o\in\O''$ then $\o=(RTR)(T\o)$.  
\end{itemize}
\end{thm}
{\bf Proof.} The measure-preserving property (i) is essentially equivalent to the output theorem 
for the stationary $M/M/1$ queue, as discussed in section 1, which follows easily from the fact 
that the Markov chain $q=s-x$ is stationary and reversible.  Property (ii) is well-known.
Properties (iii) and (iv) follow from (i) and the fact that, for $\o\in\O''$, 
$s_n=\min_{l\ge n} y_l$, $\forall n$. \hfill $\Box$

An immediate consequence of (iv) is that there exists $\O^*\subset\O$ with
$\m(\O^*)=1$ and on which $T^k$ is defined for all $k\in\Z$.
Define a mapping $\p:\O^*\to\N^\Z$ by putting $(\p\o)_k=s_0(T^k\o)$ for each $k\in\Z$.  
Denote the shift operator on $\N^\Z$ by $\theta$ and
let $\g$ be the $\theta$-invariant product measure on $\N^\Z$ with
$$\g\{\a\in\N^\Z: \a_0=x\}=(1-\rho)\rho^x\qquad x\ge 0.$$
\begin{thm}\label{main}
\ \
\begin{itemize}
\item[i.]  $\m\circ \p^{-1}=\g$.
\item[ii.]  Almost every $\o\in\O^*$ is uniquely determined by $\p\o$.
\item[iii.]  $T=\p^{-1}\circ\theta\circ\p$ almost everywhere.
\end{itemize}
\end{thm}
{\bf Proof.}  Claim (i) follows from Theorem~\ref{burke} (iii).
To prove (ii) we first note that $\o_0=(-1)^N$ where $N=\min\{k\ge 0:\ s_0(T^k\o)=0\}$.
Indeed, if $s_0(T^k\o)>0$, then $(T^{k+1}\o)_0=-(T^k\o)_0$ whereas, if $s_0(T^k\o)=0$,
then $(T^k\o)_0=1$.
By the same reasoning, for any $k\ge 0$, we have $(T^k\o)_0=(-1)^{N_k}$, 
where $N_k=\min\{l\ge 0:\ s_0(T^{k+l}\o)=0\}.$  Thus, we can recover $((T^k\o)_0,k\in\Z)$
from $\p\o$.  In exactly the same way, for any $n\in\Z$,
we can recover $((T^k\o)_n,k\in\Z)$ from the sequence $(q_n(T^k\o),k\in\Z)$,
where $q=s-x$.  Combining this observation with the identity 
$$q_{n-1}(T^k\o)=\max\{ q_n(T^k\o)+(T^{k+1}\o)_n,0\}$$ 
we see that, for any $n\le 0$ we can recover $((T^k\o)_n,k\in\Z)$ 
from $\p\o$.  In particular, we recover $(\o_n,n\le 0)$, from $\p\o$.  
A similar argument works in the other direction, starting with the observation that,
if $s_0(T^k\o)>0$, then $(T^{k+1}\o)_1=-(T^k\o)_1$ whereas, if $s_0(T^k\o)=0$,
then $(T^{k+1}\o)_1=1$;  this leads to the conclusion that $\{\o_n,n\ge 1\}$
can be recovered from $\p\o$, which completes the proof of (ii), and (iii) follows.
 \hfill $\Box$

\section{Brownian version}

Let $(X(t),t\in\R)$ be a two-sided standard Brownian motion with drift $\nu>0$
and with $X(0)=0$.  For $t\in\R$, set
$$Y(t)=2M(t)-X(t)-2M(0),\qquad M(t)=\sup_{-\infty<s\le t} X(s).$$
The continuous analogue of Burke's theorem (see, for example,~\cite{oy02b}
and references therein) states that $Y$ has the
same law as $X$ and, moreover, that $(Y(t),t\le 0)$ is independent of $M(0)$,
which is exponentially distributed with parameter $2\nu$.
It follows that the measure-preserving transformation $T$,
which maps $X$ to $Y$, has a factor which is metrically isomorphic to
the shift operator on $\R_+^\Z$, equipped with the product measure
$\varepsilon^{\otimes\Z}$, where $\varepsilon$ is the exponential
distribution on $\R_+$ with parameter $2\nu$.  However, it is not clear
in this setting whether or not $X$ can be recovered from the sequence 
$(\sup_{-\infty<s\le 0} (T^k X)(s),\ k\in\Z)$, so we cannot conclude that
$T$ is Bernoulli.  The recovery map for the discrete case, defined in the 
proof of Theorem~\ref{main}, does not have an obvious continuous 
analogue.  It is thus an open problem to determine
whether or not this transformation is Bernoulli, or even ergodic.
This is reminiscent of a (still open) question, originally posed by Marc Yor,
in relation to the following transformation of Brownian motion.
Let $(B_t,t\ge 0)$ be a standard one-dimensional Brownian motion.
It is a classical result, due to Paul L\'evy, that the process
$$\left( |B_t|-L_t^0(|B|),\ t\ge 0\right)$$
is also a standard Brownian motion, where $L_t^0(|B|)$ denotes
the local time at zero of $|B|$ up to time $t$.  Is this an ergodic
transformation?  Dubins and Smorodinsky~\cite{ds}
proved that there is a discrete version which is metrically isomorphic 
to a (one-sided) Bernoulli shift.

\section{Iterated random functions}

The classical output theorem for the $M/M/1$ queue extends quite naturally
to the more general setting of iterated random functions.  
Loosely following~\cite{df}, let $S$ be a topological space equipped with its Borel
$\sigma$-algebra, $\{ f_\t,\ \t\in\T\}$ a family of continuous functions that map $S$
to itself and $\mu$ a probability distribution on $\T$.  
Let $(\t_n,n\in\Z)$ be a sequence of random variables with common law $\kappa$.
Consider the markov chain $x=(x_n,n\ge 0)$ with state space $S$ defined
by $x_0=s$ and
\begin{equation}\label{rec}
x_n=f_{\t_{n}}(x_{n-1})=(f_{\t_n}\circ\cdots\circ f_{\t_1})(s) , \qquad n>0.
\end{equation}
We will assume that this Markov chain has reversible transition probabilities
with respect to a unique invariant probability measure.  
Now consider the {\em backward} iterations:
$$u_m=(f_{\t_1}\circ\cdots\circ f_{\t_m})(s).$$
Under certain regularity conditions, as discussed in~\cite{df}, the sequence $u_m$
converges almost surely, as $m\to\infty$, to a limiting random variable $u_\infty$
which does not depend on $s$ and which realises the invariant distribution of $x$.
We will assume that this property holds.  It follows that, for each $n\in\Z$, the limit
\begin{equation}\label{Zdef}
z_n=\lim_{m\to\infty} (f_{\t_{1+n}}\circ\cdots\circ f_{\t_{m}})(s)
\end{equation}
exists almost surely and does not depend on $s$.  
By continuity, these random variables satisfy
\begin{equation}\label{brec}
z_n=f_{\t_{n+1}}(z_{n+1}), \qquad n\in\Z,
\end{equation}
from which it follows, recalling that $x$ is has reversible transition probabilities, 
that the sequence $z=(z_n,\ n\in\Z)$ is a two-sided stationary version of $x$.  
Now suppose that, for each $s \in S$, the map $\t\mapsto (s,f_\t(s))$ is injective,
and define $F(r,s)=\t$ whenever $s=f_\t(r)$.  Then we can write
\begin{equation}\label{z}
\t_{n}=F(z_{n},z_{n-1}), \qquad n\in\Z.
\end{equation}
Define a sequence of random variables $\eta=\{\eta_n,n\in\Z\}$ by setting 
\begin{equation}\label{eta}
\eta_{n}=F(z_{n-1},z_{n}), \qquad n\in\Z,
\end{equation}  
so that 
\begin{equation}\label{brec2}
z_{n}=f_{\eta_{n}}(z_{n-1}), \qquad n\in\Z .
\end{equation}
Reversibility ensures that $\eta$ is well-defined.
\begin{thm}\label{burke2}  In the above context, $\eta$ has the same distribution as $\t$ 
and the sequence $\eta_1,\eta_2,\ldots$ is independent of $z_0$.
\end{thm}
{\bf Proof.}
The first claim follows from (\ref{z}) and (\ref{eta}),
and the fact that $z$ is stationary and reversible.
By (\ref{brec2}) we can write, almost surely, 
$$z_0=f_{\eta_0} ( f_{\eta_{-1}} (f_{\eta_{-2}}( \cdots $$ which is 
independent of $\eta_1,\eta_2,\ldots$ as required.
\hfill $\Box$

This defines a measure-preserving transformation (mapping $\theta$ to $\eta$) 
which has a Bernoulli factor.  When is it Bernoulli?
The $M/M/1$ queue corresponds to the special case where
$\T=\{-1,1\}$, $1-\kappa\{-1\}=\kappa\{1\}=q<1/2$ and 
$f_{\theta}(x)=\max\{x+\theta,0\}$.  Examples of iterated random
functions where Theorem~\ref{burke2} applies can be found in~\cite{dyson}
and~\cite{oc03}.  Further examples which arise from taking products
of random matrices, and for which the invariant measure is known explicitly,
are discussed in the paper~\cite{mtw}; 
note however that not all of these are reversible.

\noindent {\em Acknowledgements.} This research was supported by Science 
Foundation Ireland, grant number SFI 04-RP1-I512.

\end{document}